\PassOptionsToPackage{inline}{enumitem}

\documentclass[submission]{Arxiv}

\usepackage[english]{babel}
\usepackage{amsmath,amsfonts,amssymb,amsthm}
\usepackage{mathtools}
\usepackage[T1]{fontenc}
\usepackage{dsfont}
\usepackage{multicol}
\usepackage{stmaryrd}
\usepackage{mleftright}
\usepackage{csquotes}
\usepackage{tikz}
\usetikzlibrary{arrows.meta}
\usepackage{shuffle}

\setenumerate{noitemsep,topsep=0pt}

\usepackage{xcolor}
\definecolor{ColA}{HTML}{0061b7} % Dark blue.
\definecolor{ColB}{HTML}{b75600} % Dark orange.
\definecolor{FG}{RGB}{34,139,34} % Forest Green
\newcommand{\ColMix}[2]{\textcolor{ColA!#1!ColB}{#2}}
\newcommand{\ColA}[1]{\ColMix{100}{#1}}
\newcommand{\ColB}[1]{\ColMix{0}{#1}}

\setenumerate{noitemsep, topsep=6pt, itemsep=6pt, leftmargin=16pt}
\setlist[itemize]{noitemsep, topsep=6pt, itemsep=6pt, leftmargin=16pt}

\newtheorem{Theorem}{Theorem}[section]
\newtheorem{Proposition}[Theorem]{Proposition}
\newtheorem{Lemma}[Theorem]{Lemma}

\numberwithin{equation}{section}

\renewcommand{\leq}{\leqslant}
\renewcommand{\geq}{\geqslant}

\title{Two associative operads of packed words}
\author{%
    Samuele Giraudo%
    \thanks{
        \href{mailto:giraudo.samuele@uqam.ca}{\tt giraudo.samuele@uqam.ca}.
        This research has been partially supported by the projects CARPLO (ANR-20-CE40-0007)
        and ALCOHOL (ANR-19-CE40-0006) of the Agence nationale de la recherche.}%
    \addressmark{1}
    \and
    Yannic Vargas%
    \thanks{
        \href{mailto:yvargaslozada@tugraz.at}{\tt yvargaslozada@tugraz.at}.
        This research has been supported by Austrian Science Fund FWF, grant I 5788 PAGCAP.}%
    \addressmark{2}
    }
\address{%
    \addressmark{1} Université du Québec à Montréal, LACIM, Pavillon Président-Kennedy, 201
    Avenue du Président-Kennedy, Montréal, H2X~3Y7, Canada. \\
    \addressmark{2} Institute of Discrete Mathematics, Technische Universität Graz, 
    Austria, Graz, Austria.}
\abstract{%
    The associative operad is a central structure in operad theory, defined on the linear
    span of the set of permutations. We build two analogs of the associative operad on the
    linear span of the set of packed words which turn out to be set-theoretical. By seeing a
    packed word as a surjective map between two finite sets, our first operad is graded by
    the cardinality of the domain and the second one, by the cardinality of the codomain. In
    the same way as the associative operad of permutations contains as quotients the
    duplicial and interstice operads, we derive similar structures for our operads of packed
    words. We propose also an analogue of Dynkin idempotent of Zie algebras in this context
    of operads of packed words.
}
\received{\today}

\usepackage[backend=bibtex]{biblatex}
\newcommand{\st}{\operatorname{st}}

\newcommand{\Def}[1]{\ColMix{50}{\em #1}}
\newcommand{\Par}[1]{\mleft(#1\mright)}
\newcommand{\OEIS}[1]{\href{http://oeis.org/#1}{{\bf #1}}}

\newcommand{\N}{\mathbb{N}}

\newcommand{\PositiveN}{\mathbb{P}}
\newcommand{\Length}{\ell}
\newcommand{\SymmetricAction}{\cdot}
\newcommand{\InversionSet}{\text{Inv}}
\newcommand{\Inversions}{\text{inv}}

\newcommand{\Increment}[2]{\uparrow_{#1}^{#2}}
\newcommand{\SetPermutations}{\mathfrak{S}}
\newcommand{\SetPackedWords}{\mathfrak{P}}
\newcommand{\PW}{\mathsf{PW}}
\newcommand{\As}{\mathsf{As}}
\newcommand{\Dup}{\mathsf{Dup}}
\newcommand{\Interstice}{\mathsf{I}}
\newcommand{\RightPAs}{\mathsf{{PAs}^{\rightarrow}}}
\newcommand{\LeftPAs}{\mathsf{{PAs}^{\leftarrow}}}
\newcommand{\RightAs}{\As^{\rightarrow}}
\newcommand{\LeftAs}{\As^{\leftarrow}}
\newcommand{\PredicateSylv}{P_{\mathrm{S}}}
\newcommand{\PredicateCoSylv}{\bar{P_{\mathrm{S}}}}
\newcommand{\PredicateComm}{P_{\mathrm{C}}}
\newcommand{\PredicateBaxter}{P_{\mathrm{B}}}
\newcommand{\PredicateHypoplactic}{P_{\mathrm{H}}}
\newcommand{\PredicatePlactic}{P_{\mathrm{P}}}
\newcommand{\Equiv}{\equiv}
\newcommand{\First}[1]{\mathrm{first}_{#1}}
\newcommand{\AdjacencyRelation}{\leftrightarrow}
\newcommand{\NCSF}{\mathbf{NCSF}}
\newcommand{\Conc}{.}
\newcommand{\Tree}{\mathfrak{t}}
\newcommand{\TreeS}{\mathfrak{s}}
\newcommand{\TreeR}{\mathfrak{r}}
\newcommand{\BasisE}{\mathsf{E}}
\newcommand{\BasisF}{\mathsf{F}}
\newcommand{\BasisM}{\mathsf{M}}
\newcommand{\record}{\operatorname{\mathrm{rec}}}
\newcommand{\WQSym}{\mathsf{WQSym}}
\newcommand{\FQSym}{\mathsf{FQSym}}

\tikzstyle{Centering}=[{baseline={([yshift=-0.5ex]current bounding box.center)}}]
\tikzstyle{Node}=[circle,draw=ColA!80,fill=ColA!8,inner sep=1pt, minimum size=2mm, thick,
    font=\scriptsize]
\tikzstyle{Edge}=[draw=ColB!80,cap=round,thick,rounded corners=2.5pt]
\tikzstyle{Leaf}=[rectangle,draw=ColA!70,fill=ColA!16,inner sep=0pt,minimum size=1mm,thick]

\begin{document}
\maketitle

%%%%%%%%%%%%%%%%%%%%%%%%%%%%%%%%%%%%%%%%%%%%%%%%%%%%%%%%%%%%%%%%%%%%%%%%%%%%%%%%%%%%%%%%%%%%
%%%%%%%%%%%%%%%%%%%%%%%%%%%%%%%%%%%%%%%%%%%%%%%%%%%%%%%%%%%%%%%%%%%%%%%%%%%%%%%%%%%%%%%%%%%%
%%%%%%%%%%%%%%%%%%%%%%%%%%%%%%%%%%%%%%%%%%%%%%%%%%%%%%%%%%%%%%%%%%%%%%%%%%%%%%%%%%%%%%%%%%%%
\section*{Introduction}
The associative operad $\As$ is an important algebraic structure which plays a central role
in the theory of operads. This operad is defined on the linear span of the set of
permutations and its partial composition consists of inserting a permutation into another
one, interpreted as permutation matrices. The first reason justifying the importance of
$\As$ is that it intervenes in a crucial way in the description of the axioms of symmetric
operads. A second reason, shown by Aguiar and Livernet~\cite{AL2007}, relates to its
richness from a combinatorial point of view. Indeed, $\As$ admits a basis defined using the
left weak order on permutations which possesses the nice property that the partial
composition of two elements of this basis is a sum of an interval of this partial order.
Furthermore, This operad enjoys many interesting properties since it contains not only the
Lie operad as a suboperad, but also, as quotients, the duplicial operad of binary
trees~\cite{BF03} and the interstice operad on two generators of binary
words~\cite{Chapoton2002,CG22}.
\medbreak

These three operads on permutations, binary trees, and binary words form a hierarchy very
similar to a well-known hierarchy of combinatorial Hopf algebras involving the same spaces
of combinatorial objects, namely, the Malvenuto-Reutenauer Hopf algebra of
permutations~\cite{MR1995,DHT02}, the Loday-Ronco Hopf algebra of binary
trees~\cite{LR1998}, and the noncommutative symmetric functions Hopf
algebra~\cite{GKLLRT94}. Interestingly enough, there is a generalization of the
Malvenuto-Reutenauer Hopf algebra on the linear span of packed words.  This  Hopf algebra
has been introduced by Hivert~\cite{Hiv99} when he considered a notion of word
quasi-symmetric functions. This construction is natural since while permutations are
bijections, packed words are surjections. In this context, the analog of the Loday-Ronco
Hopf algebra involves Schröder trees~\cite{NT20} and the analog of the noncommutative
symmetric functions  Hopf algebra involves ternary words~\cite{NT20}. The starting point of
the present work is to explore whether such a natural generalization of $\As$ exists and if
it leads to a similar hierarchy of operads.
\medbreak

Our main contribution consists of the introduction of two different generalizations of $\As$
on the linear span of packed words. By seeing a packed word as its matrix, we obtain a right
version $\RightPAs$ consisting of inserting a packed word matrix at a given position into
another one, and a left version $\LeftPAs$ consisting of inserting several copies of a
packed word matrix onto given values. These operads also differ in the way the arity of a
packed word is defined: in $\RightPAs$ (resp.\ $\LeftPAs$), the arity of a packed word is
the cardinality of its domain (resp.\ codomain). As they are not isomorphic as graded
spaces, $\RightPAs$ and $\LeftPAs$ are not isomorphic as operads. Moreover, $\LeftPAs$ is a
symmetric operad while $\RightPAs$ is not. Besides, the operad $\LeftPAs$ is not
combinatorial in the sense that it admits infinitely many elements of any given arity $n
\geq 1$. Nevertheless, despite this fact, this operad is rich from a combinatorial point of
view since it admits several quotients using well-known equivalence relations on packed
words, like the sylvester~\cite{HNT05}, hypoplactic~\cite{KT97}, or Baxter~\cite{Gir12}
congruences.
\medbreak

This work is presented as follows. Section~\ref{sec:preliminaries} contains fundamental
notions about the main combinatorial objects and operads appearing in this work. In
Section~\ref{sec:operads_packed_words}, we introduce and construct the operads $\RightPAs$
and $\LeftPAs$ and present their first properties. Section~\ref{sec:quotients} is devoted to
the study of some of the quotients of these two operads of packed words which involve other
families of combinatorial objects. Section~\ref{sec:dynkin} contains the construction of an
analogue of the classical Dynkin idempotents for Zie algebras, introduced in \cite{AM2017}.
Finally, Section~\ref{sec:open_questions} contains several open questions as well as avenues
for future research in continuation of this work.
\medbreak

%%%%%%%%%%%%%%%%%%%%%%%%%%%%%%%%%%%%%%%%%%%%%%%%%%%%%%%%%%%%%%%%%%%%%%%%%%%%%%%%%%%%%%%%%%%%
{\it General notations and conventions.}
For any integer $i$, $[i]$ denotes the set $\{1, \dots, i\}$. For any set $A$, $A^*$ is the
set of words on $A$. For any $w \in A^*$, $\Length(w)$ is the length of $w$, and for any $i
\in [\Length(w)]$, $w(i)$ is the $i$-th letter of $w$. The only word of length $0$ is the
empty word $\epsilon$. For any $1 \leq i \leq j \leq \Length(w)$, $w(i, j)$ is the word
$w(i) \dots w(j)$. Given two words $w$ and $w'$, the  concatenation of $w$ and $w'$ is
denoted by $ww'$ or by~$w \Conc w'$.
\medbreak

%%%%%%%%%%%%%%%%%%%%%%%%%%%%%%%%%%%%%%%%%%%%%%%%%%%%%%%%%%%%%%%%%%%%%%%%%%%%%%%%%%%%%%%%%%%%
%%%%%%%%%%%%%%%%%%%%%%%%%%%%%%%%%%%%%%%%%%%%%%%%%%%%%%%%%%%%%%%%%%%%%%%%%%%%%%%%%%%%%%%%%%%%
%%%%%%%%%%%%%%%%%%%%%%%%%%%%%%%%%%%%%%%%%%%%%%%%%%%%%%%%%%%%%%%%%%%%%%%%%%%%%%%%%%%%%%%%%%%%
\section{Preliminaries} \label{sec:preliminaries}

%%%%%%%%%%%%%%%%%%%%%%%%%%%%%%%%%%%%%%%%%%%%%%%%%%%%%%%%%%%%%%%%%%%%%%%%%%%%%%%%%%%%%%%%%%%%
%%%%%%%%%%%%%%%%%%%%%%%%%%%%%%%%%%%%%%%%%%%%%%%%%%%%%%%%%%%%%%%%%%%%%%%%%%%%%%%%%%%%%%%%%%%%
\subsection{Packed words and related objects}
Let $\PositiveN$ be the set $\N \setminus \{0\}$. For any $u \in \PositiveN^*$ and $\alpha,
\beta \in \N$, let $\Increment{\alpha}{\beta}(u)$ be the word obtained by incrementing by
$\alpha$ the letters of $u$ which are greater than $\beta$. For instance,
$\Increment{2}{4}(124\ColB{5}4\ColB{6}) = 124\ColB{7}4\ColB{8}$. Given $w \in \PositiveN^*$,
the \Def{alphabet} of $w$ is the set $\text{Alph}(w)$ formed by all different letters
appearing in $w$. An \Def{inversion} of a word $w \in \PositiveN^*$ is a pair $(i,j)$ such
that $1 \leq i < j \leq \Length(w)$ and $w(i)>w(j)$. The \Def{set of inversions} of $w$ is
denoted by $\InversionSet(w)$ and the \Def{number of inversions} of $w$ is
$\Inversions(w):=|\InversionSet(w)|$.
\medbreak

The \Def{standardization map} $\st: \PositiveN^* \to \mathfrak{S}$ is the function sending a
word $w$ of length $n$ to a permutation $\st(w) \in \mathfrak{S}_n$,  obtained by
iteratively scanning $w$ from left to right, and labelling $1,2,3, \hdots$ the occurrences
of its smallest letter, then numbering the occurrences of the next one, and so on. For
instance, if $x<y<z$ are letters in $\PositiveN^*$, then $\text{st}(yyxzzxzyx)=451782963 \in
\mathfrak{S}_9$. The standardization of a word $w$ is the unique permutation of
$\SetPermutations_n$ preserving the relative order of $w$; that is, $w$ and $\st(w)$ have
the same inversions (see \cite{Hivert2006}).
\medbreak

A \Def{packed word} is a word $w \in \PositiveN^*$ satisfying $i - 1 \in \text{Alph}(w)$
whenever $i \in \text{Alph}(w)$, for every $i \geq 2$. This implies that the alphabet of a
packed word $w$ is such that $\text{Alph}(w)=[k]$, for some $k \geq 0$. Also, if $n$ is the
length of $w$, the definition of packed words necessarily implies that $k \leq n$. For every
$k,n \geq 0$, let $\mathfrak{P}_n[k]$ be the set of all packed words of alphabet $[k]$ and
length $n$. For every $n \geq 0$, $\SetPermutations_n:=\mathfrak{P}_n[n]$ is the set of
\Def{permutations} of the set $[n]$.
\medbreak

A \Def{Schröder tree} is a planar rooted tree where each internal node has at least two
children. For every $k,n \geq 0$, let $\mathfrak{T}_n[k]$ be the set of all Schröder trees
with $n+1$ leaves and $k$ internal nodes. Observe that $\mathfrak{B}_n := \mathfrak{T}_n[n]$
is the set of \Def{binary trees} with $n+1$ leaves (or equivalently, binary trees with $n$
internal nodes), for every $n \geq 0$.
\medbreak

Given $k,n \geq 0$ and $\mathfrak{X}_k[n] \in \{ \mathfrak{P}_k[n], \mathfrak{T}_k[n]\}$, we
put
\begin{equation}
    \mathfrak{X}_n:=\bigcup_{0 \leq k \leq n} \mathfrak{X}_n[k],
    \quad
    \mathfrak{X}[k]:=\bigcup_{n \geq k} \mathfrak{X}_n[k],
    \quad \text{and} \quad
    \mathfrak{X}:= \bigcup_{n \geq 0}\mathfrak{X}_n.
\end{equation}
Here, $\mathfrak{X}_n$ and $\mathfrak{X}[k]$  correspond to the sets of our given
combinatorial objects classified by \Def{length} and \Def{alphabet}, whenever $\mathfrak{X}
\in \{\mathfrak{P}, \mathfrak{T}\}$ represent packed words or Schröder trees, respectively.
\medbreak

Recall that the symmetric group acts on the vector space of non-commutative polynomials
$\mathbb{K}\langle \mathsf{X}\rangle$, for any countable set $\mathsf{X}$, where
$\mathbb{K}$ is any field of characteristic zero. When $\mathsf{X}=\PositiveN$, this action
restricts to packed words: if $w \in \mathfrak{P}_n[k]$ and $\sigma \in \mathfrak{S}_n$,
then
\begin{equation} \label{equ:right_action_permutation}
    w \cdot \sigma:= w(\sigma(1))\cdots w(\sigma(n)) \in \mathfrak{P}_n[k].
\end{equation}
Indeed, if $w$ is a packed word, any permutation of its letters leads to a word satisfying
the condition on the definition of packed words.
\medbreak

For every packed word $w$, let $\chi(w)$ be the unique weakly increasing word obtained after
rearrangement of $w$. We call $\chi(w)$ the \Def{composition type} of $w$. In particular,
$\chi(w)$ is again a packed word of same alphabet and length than $w$. Every weakly
increasing packed word $w \in \mathfrak{P}_n[k]$ can be encoded by the integer composition
$\mathsf{c} := c_1c_2 \dots c_k$ of $n$ such that every $i$ appears $c_i$ times in $w$.
Equivalently, through the classical bijection between subsets of $[n-1]$ and binary strings
of length $n-1$, the composition type $\chi(w)$ of a packed word of length $n$ can be also
described as a word of length $n-1$ on the alphabet $\{1,2\}$. For instance, if $u =
2311223$ and $v=221$, then
\begin{equation}\label{Example_chi}
    \chi(u) = 1122233 \leftrightarrow \ColA{-}\ColB{+}\ColA{--}\ColB{+}\ColA{-}
    \leftrightarrow \ColA{1}\ColB{2}\ColA{11}\ColB{2}\ColA{1}
    \enspace \text{and} \enspace
    \chi(v) = 122  \leftrightarrow \ColB{+}\ColA{-} \leftrightarrow \ColB{2}\ColA{1}.
\end{equation}
\medbreak

The following lemma, due to Hivert, relates a packed word with its composition type.
\medbreak

\begin{Lemma}[Hivert, \cite{Hivert2006}]\label{stHiv2}
     For any $w \in \mathfrak{P}_n[k]$, 
     \begin{math}
         w = \chi(w)\cdot \text{st}(w).
     \end{math}
    Moreover, $\text{st}(w)$ is the smallest permutation in $\mathfrak{S}_n$ for the right
    weak order satisfying the above property.
\end{Lemma}
\medbreak

This result implies that every packed word $w$ is encoded by its composition type and its
standardization: the composition type tells us how many times each letter appears in $w$,
while $\text{st}(w)$ encodes the relative order of the appearance of the letters in~$w$.
\medbreak

%%%%%%%%%%%%%%%%%%%%%%%%%%%%%%%%%%%%%%%%%%%%%%%%%%%%%%%%%%%%%%%%%%%%%%%%%%%%%%%%%%%%%%%%%%%%
%%%%%%%%%%%%%%%%%%%%%%%%%%%%%%%%%%%%%%%%%%%%%%%%%%%%%%%%%%%%%%%%%%%%%%%%%%%%%%%%%%%%%%%%%%%%
\subsection{Associative operad and quotients}
We follow the usual notations about symmetric unital operads~\cite{LV2012} (called simply
\Def{operads} here). Here are four important examples of (nonsymmetric) operads.
\begin{itemize}
    \item The \Def{right associative operad} is the symmetric operad $\RightAs$ (also
    written $\As$) wherein for any $n \geq 1$, $\RightAs(n)$ is the linear span of the set
    $\mathfrak{S}_n$. The set $\{\BasisE_\sigma : \sigma \in \mathfrak{S}_n, \, n \geq 1\}$
    is a basis of $\RightAs$. Given $\alpha \in \mathfrak{S}_n$, $\beta \in \mathfrak{S}_m$
    and $i \in [n]$, let 
    \begin{equation}\label{defB}
        B_i(\alpha, \beta) :=
            \Increment{m - 1}{\alpha(i)}(\alpha(1, i - 1))
            \; \Conc \;
            \Increment{\alpha(i) - 1}{0}(\beta)
            \; \Conc \;
            \Increment{m - 1}{\alpha(i) - 1}(\alpha(i + 1, n)).
    \end{equation}
    For instance, 
    \begin{math}
        B_5(\ColA{3612} \ColB{4} \ColA{57}, \ColB{231}) = \ColA{3812} \ColB{564} \ColA{79}.
    \end{math}
    The \ColB{red} labels reflect the \ColB{red} permutation $231$ inserted into the
    \ColA{blue} permutation $3612457$ onto the letter $\ColB{4}$ at the $5$-th position. The
    permutation $B_i(\alpha, \beta)$ is sometimes called the \Def{block permutation}
    associated to $\alpha$ and $\beta$ (see~\cite{LV2012}). The partial composition of
    $\RightAs$ satisfies $\BasisE_\sigma \circ_i \BasisE_\nu = \BasisE_{B_i(\sigma, \nu)}$
    and the action of the symmetric group $\mathfrak{S}_n$ satisfies $\BasisE_\sigma \cdot
    \nu = \BasisE_{\sigma \circ \nu}$ where $\circ$ is the composition of permutations.
    \item The \Def{left associative operad} is the symmetric operad $\LeftAs$ defined on the
    same space as the one of $\RightAs$, seen on the same $\BasisE$-basis. The partial
    composition of $\LeftPAs$ satisfies $\BasisE_\sigma \circ_i \BasisE_\nu = \BasisE_\pi$
    where $\pi$ is the permutation obtained by replacing in $\Increment{m - 1}{i}(\sigma)$
    the letter $i$ by $\Increment{i - 1}{0}(\nu)$, where $m$ is the greatest letter of
    $\nu$. The action of the symmetric group $\mathfrak{S}_n$ satisfies $\BasisE_\sigma
    \cdot \nu = \BasisE_{\nu^{-1} \circ \sigma}$ where $\circ$ is the composition of
    permutations. The operads $\RightAs$ and $\LeftAs$ are isomorphic through the linear map
    $\phi : \RightAs \to \LeftAs$ satisfying $\phi\left(\BasisE_\sigma\right) =
    \BasisE_{\sigma^{-1}}$.
    \item The \Def{duplicial operad}~\cite{BF03} is the nonsymmetric operad $\Dup$ wherein
    for any $n \geq 1$, $\Dup(n)$ is the linear span of the set $\mathfrak{B}_n$. The set
    $\{\BasisE_\Tree : \Tree \in \mathfrak{B}_n, \, n \geq 1\}$ is a basis of $\Dup$. The
    partial composition of $\Dup$ satisfies $\BasisE_\Tree \circ_i \BasisE_\TreeS =
    \BasisE_\TreeR$ where $\TreeR$ is the binary tree obtained by replacing the $i$-th
    internal node $u$ of $\Tree$ for the infix traversal by a copy of $\TreeS$ and by
    grafting the left subtree of $u$ on the leftmost leaf of the copy and by grafting the
    right subtree of $u$ on the rightmost least of the copy. For instance, in $\Dup$, we
    have
    \begin{equation}
        \BasisE_{
            \scalebox{.85}{
            \begin{tikzpicture}[xscale=.16,yscale=.12,Centering]
                \node[Leaf](0)at(0.00,-9.00){};
                \node[Leaf](10)at(10.00,-12.00){};
                \node[Leaf](12)at(12.00,-9.00){};
                \node[Leaf](14)at(14.00,-9.00){};
                \node[Leaf](2)at(2.00,-9.00){};
                \node[Leaf](4)at(4.00,-6.00){};
                \node[Leaf](6)at(6.00,-9.00){};
                \node[Leaf](8)at(8.00,-12.00){};
                \node[Node](1)at(1.00,-6.00){};
                \node[Node,draw=ColB,fill=ColB!50](11)at(11.00,-3.00){};
                \node[Node](13)at(13.00,-6.00){};
                \node[Node](3)at(3.00,-3.00){};
                \node[Node](5)at(5.00,0.00){};
                \node[Node](7)at(7.00,-6.00){};
                \node[Node](9)at(9.00,-9.00){};
                \draw[Edge](0)--(1);
                \draw[Edge](1)--(3);
                \draw[Edge](10)--(9);
                \draw[Edge](11)--(5);
                \draw[Edge](12)--(13);
                \draw[Edge](13)--(11);
                \draw[Edge](14)--(13);
                \draw[Edge](2)--(1);
                \draw[Edge](3)--(5);
                \draw[Edge](4)--(3);
                \draw[Edge](6)--(7);
                \draw[Edge](7)--(11);
                \draw[Edge](8)--(9);
                \draw[Edge](9)--(7);
                \node(r)at(5.00,3.5){};
                \draw[Edge](r)--(5);
            \end{tikzpicture}}}
        \enspace \circ_{\ColB{6}} \enspace
        \BasisE_{
            \scalebox{.85}{
            \begin{tikzpicture}[xscale=.16,yscale=.15,Centering]
                \node[Leaf](0)at(0.00,-4.50){};
                \node[Leaf](2)at(2.00,-4.50){};
                \node[Leaf](4)at(4.00,-4.50){};
                \node[Leaf](6)at(6.00,-6.75){};
                \node[Leaf](8)at(8.00,-6.75){};
                \node[Node,draw=ColB,fill=ColB!50](1)at(1.00,-2.25){};
                \node[Node,draw=ColB,fill=ColB!50](3)at(3.00,0.00){};
                \node[Node,draw=ColB,fill=ColB!50](5)at(5.00,-2.25){};
                \node[Node,draw=ColB,fill=ColB!50](7)at(7.00,-4.50){};
                \draw[Edge](0)--(1);
                \draw[Edge](1)--(3);
                \draw[Edge](2)--(1);
                \draw[Edge](4)--(5);
                \draw[Edge](5)--(3);
                \draw[Edge](6)--(7);
                \draw[Edge](7)--(5);
                \draw[Edge](8)--(7);
                \node(r)at(3.00,3){};
                \draw[Edge](r)--(3);
            \end{tikzpicture}}}
        \enspace = \enspace
        \BasisE_{
            \scalebox{.85}{
            \begin{tikzpicture}[xscale=.15,yscale=.100,Centering]
                \node[Leaf](0)at(0.00,-10.50){};
                \node[Leaf](10)at(10.00,-17.50){};
                \node[Leaf](12)at(12.00,-10.50){};
                \node[Leaf](14)at(14.00,-10.50){};
                \node[Leaf](16)at(16.00,-14.00){};
                \node[Leaf](18)at(18.00,-17.50){};
                \node[Leaf](2)at(2.00,-10.50){};
                \node[Leaf](20)at(20.00,-17.50){};
                \node[Leaf](4)at(4.00,-7.00){};
                \node[Leaf](6)at(6.00,-14.00){};
                \node[Leaf](8)at(8.00,-17.50){};
                \node[Node](1)at(1.00,-7.00){};
                \node[Node](3)at(3.00,-3.50){};
                \node[Node](5)at(5.00,0.00){};
                \node[Node](7)at(7.00,-10.50){};
                \node[Node](9)at(9.00,-14.00){};
                \node[Node,draw=ColB,fill=ColB!50](11)at(11.00,-7.00){};
                \node[Node,draw=ColB,fill=ColB!50](13)at(13.00,-3.50){};
                \node[Node,draw=ColB,fill=ColB!50](15)at(15.00,-7.00){};
                \node[Node,draw=ColB,fill=ColB!50](17)at(17.00,-10.50){};
                \node[Node](19)at(19.00,-14.00){};
                \draw[Edge](0)--(1);
                \draw[Edge](1)--(3);
                \draw[Edge](10)--(9);
                \draw[Edge](11)--(13);
                \draw[Edge](12)--(11);
                \draw[Edge](13)--(5);
                \draw[Edge](14)--(15);
                \draw[Edge](15)--(13);
                \draw[Edge](16)--(17);
                \draw[Edge](17)--(15);
                \draw[Edge](18)--(19);
                \draw[Edge](19)--(17);
                \draw[Edge](2)--(1);
                \draw[Edge](20)--(19);
                \draw[Edge](3)--(5);
                \draw[Edge](4)--(3);
                \draw[Edge](6)--(7);
                \draw[Edge](7)--(11);
                \draw[Edge](8)--(9);
                \draw[Edge](9)--(7);
                \node(r)at(5.00,4){};
                \draw[Edge](r)--(5);
        \end{tikzpicture}}}.
    \end{equation}
    \item For any $s \geq 1$, the \Def{$s$-interstice operad}~\cite{Chapoton2002,CG22} is
    the nonsymmetric operad $\Interstice_s$ wherein for any $n \geq 1$, $\Interstice_s(n)$
    is the linear span of the set $[s]^n$. The set $\{\BasisE_u : u \in [s]^n, n \geq 1\}$
    is a basis of $\Interstice_s$. The partial composition of $\Interstice_s$ satisfies
    \begin{equation}\label{def_sInterstice}
        \BasisE_u \circ_i \BasisE_v
        := \BasisE_{u(1, i - 1) 
        \; \Conc \; v \; \Conc \;
        u(i, \Length(u))}.
    \end{equation}
    For instance, in $\Interstice_2$, we have
    \begin{math}
        \BasisE_{\ColA{121 121}} \circ_4 \BasisE_{\ColB{21}}
        =
        \BasisE_{\ColA{121} \; \ColB{21} \; \ColA{121}}.
    \end{math}
\end{itemize}
As shown in~\cite{AL2007}, $\Dup$ and $\Interstice_2$ are nonsymmetric operad quotients of
$\RightPAs$.
\medbreak

%%%%%%%%%%%%%%%%%%%%%%%%%%%%%%%%%%%%%%%%%%%%%%%%%%%%%%%%%%%%%%%%%%%%%%%%%%%%%%%%%%%%%%%%%%%%
%%%%%%%%%%%%%%%%%%%%%%%%%%%%%%%%%%%%%%%%%%%%%%%%%%%%%%%%%%%%%%%%%%%%%%%%%%%%%%%%%%%%%%%%%%%%
%%%%%%%%%%%%%%%%%%%%%%%%%%%%%%%%%%%%%%%%%%%%%%%%%%%%%%%%%%%%%%%%%%%%%%%%%%%%%%%%%%%%%%%%%%%%
\section{Operadic structures on packed words} \label{sec:operads_packed_words}
We introduce two generalizations $\RightPAs$ and $\LeftPAs$ respectively of $\RightAs$ and
$\LeftAs$ on the linear span of the set of packed words seen through two different
graduations.
\medbreak

%%%%%%%%%%%%%%%%%%%%%%%%%%%%%%%%%%%%%%%%%%%%%%%%%%%%%%%%%%%%%%%%%%%%%%%%%%%%%%%%%%%%%%%%%%%%
%%%%%%%%%%%%%%%%%%%%%%%%%%%%%%%%%%%%%%%%%%%%%%%%%%%%%%%%%%%%%%%%%%%%%%%%%%%%%%%%%%%%%%%%%%%%
\subsection{Right version}
Let $\RightPAs$ be the graded space such that for any $n \geq 1$, $\RightPAs(n)$ is the
linear span of the set $\SetPackedWords_n$. The set $\{\BasisE_u : u \in \SetPackedWords_n,
n \geq 1\}$ is a basis of $\RightPAs$. Let us endow $\RightPAs$ with the operations
$\circ_i$ defined, for any $u \in \SetPackedWords_n[r]$, $i \in [n]$, and $v \in
\SetPackedWords_{m}[s]$, by
\begin{equation}\label{RightPAs_Def}
    \BasisE_u \circ_i \BasisE_v :=
        \BasisE_{
        \Increment{s - 1}{u(i)}(u(1, i - 1))
        \; \Conc \;
        \Increment{u(i) - 1}{0}(v)
        \; \Conc \;
        \Increment{s - 1}{u(i) - 1}(u(i + 1, n))}.
\end{equation}
For instance,
\begin{math}
    \BasisE_{\ColA{{\bf 2} 311} \ColB{\bf 2} \ColA{ {\bf 2}3}}
    \circ_4
    \BasisE_{\ColB{221}}
    =
    \BasisE_{\ColA{{\bf 2} 411} \; \ColB{332} \; \ColA{ {\bf 3}4}}.
\end{math}
\medbreak

Intuitively, the partial composition $\BasisE_u \circ_i \BasisE_v$ of $\RightPAs$ is similar
to the one $\RightAs$ but with the difference that the occurrences of $u(i)$ in $u$ having a
position greater than $i$ are incremented by $s - 1$ where $s$ is the maximal value of $v$.
In terms of permutation matrices, we have

\begin{footnotesize}
\begin{equation}
    \BasisE_{
    \begin{bmatrix}
        \ColA{0} & \ColA{0} & \ColA{1} &  \ColA{1} & 0 & \ColA{0} & \ColA{0} \\
        \ColA{1} & \ColA{0} & \ColA{0}  & \ColA{0} & \ColB{\bf 1} & \ColA{1} & \ColA{0} \\
        \ColA{0} & \ColA{1} & \ColA{0} &  \ColA{0} & 0 & \ColA{0} & \ColA{1}
    \end{bmatrix}}
    \circ_4 
    \BasisE_{
    \begin{bmatrix}
        \ColB{0} & \ColB{0} & \ColB{1} \\
        \ColB{1} & \ColB{1} & \ColB{0}
    \end{bmatrix}}
    =
    \BasisE_{
    \begin{bmatrix}
        \ColA{0} & \ColA{0} & \ColA{1} & \ColA{1} & 0 & 0 & 0 & \ColA{0} & \ColA{0} \\
        \ColA{1} & \ColA{0} & \ColA{0} & \ColA{0} & \ColB{0} & \ColB{0} & \ColB{1} & 0 & 0 \\
        0 & 0 & 0 & 0 & \ColB{1} & \ColB{1} & \ColB{0} & \ColA{1} & \ColA{0} \\
        \ColA{0} & \ColA{1} & \ColA{0} & \ColA{0} & 0 & 0 & 0 & \ColA{0} & \ColA{1}
    \end{bmatrix}}.
\end{equation}
\end{footnotesize}
\medbreak

\begin{Theorem} \label{thm:rigth_packed_word_associative_operad}
    The graded space $\RightPAs$, endowed with the partial composition maps $\circ_i$, is a
    nonsymmetric unital operad.
\end{Theorem}
\medbreak

We call $\RightPAs$ the \Def{right associative operad of packed words}. Remark that
in~\cite{Gir2015} and~\cite{McCS2003}, other operads involving the same space of packed
words are constructed. The operads $\RightPAs$ and $\PW$, which is defined
in~\cite{Gir2015}, are not isomorphic. Indeed, a simple inspection shows that any minimal
generating set of $\PW$ contains two elements of arity $3$ while any minimal generating set
$\RightPAs$ has no element of this arity.
\medbreak

For any $u \in \mathfrak{P}_k$ and $i \in [k]$, the \Def{record} $\text{rec}_i(u)$ of $u$ at
$i$ is $\text{rec}_i(w):=\text{st}(u)(i)$. For example, taking $u=2311{\bf 2}23$ and $i=5$,
we have $\st(u)=3612{\bf 4}57$, so $\text{rec}_5(u)=\st(u)(5)=4$. In what follows, we
consider that the maps $\record$ and $\st$ are extended linearly on the graded space
$\RightPAs$ through its $\BasisE$-basis. In the same way, we extend linearly the action
$\cdot$ of~\eqref{equ:right_action_permutation} on $\RightPAs$ through its $\BasisE$-basis.
\medbreak

As shown by the next result, the operad $\RightPAs$ relates well with $\RightAs$ and
$\Interstice_2$.
\medbreak

\begin{Theorem} \label{thm:decopmposition_of_rigth_packed_word_associative_operad}
    Let $n \in \mathbb{N}$, $u \in \mathfrak{P}_n$, and $v \in \mathfrak{P}$. For any $i \in
    [n]$,
    \begin{equation}\label{Relating_B}
        \BasisE_u \circ_i \BasisE_v
        =
        \left(\BasisE_{\chi(u)} \circ_{\record_i(u)} \BasisE_{\chi(v)}\right)
        \cdot B_i(\st(u), \st(v)),
    \end{equation}
    where the partial composition $\circ_{\record_i(u)}$ is the one of $\Interstice_2$. In
    particular, $\chi(\BasisE_u \circ_i \BasisE_v) = \BasisE_{\chi(u)} \circ_{\record_i(u)}
    \BasisE_{\chi(v)}$ and $\st\left(\BasisE_u \circ_i \BasisE_v\right) =
    \BasisE_{B_i(\st(u),\, \st(v))}$.
\end{Theorem}
\medbreak

Following with our example let $u=2311{\bf 2}23$ and $v = 221$. Notice that $\st(u) =
3612457$, so that $\record_5(u)=4$. Also, $\st(v)= 231$. Therefore,
\begin{equation}
    \left(\BasisE_{\chi(u)} \circ_{\record_i(u)} \BasisE_{\chi(v)}\right)
        \cdot B_i(\st(u), \st(v))
        =
        \BasisE_{112233344} \cdot 381256479 = \BasisE_{241133234},
\end{equation}
which agrees with our previous example at the begining of this section.
\medbreak

%%%%%%%%%%%%%%%%%%%%%%%%%%%%%%%%%%%%%%%%%%%%%%%%%%%%%%%%%%%%%%%%%%%%%%%%%%%%%%%%%%%%%%%%%%%%
%%%%%%%%%%%%%%%%%%%%%%%%%%%%%%%%%%%%%%%%%%%%%%%%%%%%%%%%%%%%%%%%%%%%%%%%%%%%%%%%%%%%%%%%%%%%
\subsection{Left version}
Let $\LeftPAs$ be the graded space such that for any $n \geq 1$, $\LeftPAs(n)$ is the linear
span of the set $\SetPackedWords[n]$. The set $\{\BasisE_u : u \in \SetPackedWords[n], n
\geq 1\}$ is a basis of $\LeftPAs$. Let us endow $\LeftPAs$ with the operations $\circ_i$
defined, for any $u \in \SetPackedWords[n]$, $i \in [n]$, and $v \in \SetPackedWords[m]$, by
\begin{math}
    \BasisE_u \circ_i \BasisE_v = \BasisE_w
\end{math}
where $w$ is the word obtained by replacing in $\Increment{m - 1}{i}(u)$ all occurrences of
$i$ by $\Increment{i - 1}{0}(v)$. For instance,
\begin{equation}
    \BasisE_{\ColA{2} \ColB{3} \ColA{1} \ColB{3} \ColA{14}}
    \circ_{\ColB{3}}
    \BasisE_{\ColB{122}}
    =
    \BasisE_{\ColA{2} \; \ColB{344} \; \ColA{1} \; \ColB{344} \; \ColA{15}}.
\end{equation}
In terms of permutation matrices, we have

\begin{footnotesize}
\begin{equation}
    \BasisE_{
    \begin{bmatrix}
        \ColA{0} & 0 & \ColA{1} & 0 & \ColA{1} & \ColA{0} \\
        \ColA{1} & 0 & \ColA{0} & 0 & \ColA{0} & \ColA{0} \\
        \ColA{0} & \ColB{\bf 1} & \ColA{0} & \ColB{\bf 1} & \ColA{0} & \ColA{0} \\
        \ColA{0} & 0 & \ColA{0} & 0 & \ColA{0} & \ColA{1}
    \end{bmatrix}}
    \circ_3
    \BasisE_{
    \begin{bmatrix}
        \ColB{1} & \ColB{0} & \ColB{0} \\
        \ColB{0} & \ColB{1} & \ColB{1}
    \end{bmatrix}}
    =
    \BasisE_{
    \begin{bmatrix}
        \ColA{0} & 0 & 0 & 0 & \ColA{1} & 0 & 0 & 0 & \ColA{1} & \ColA{0} \\
        \ColA{1} & 0 & 0 & 0 & \ColA{0} & 0 & 0 & 0 & \ColA{0} & \ColA{0} \\
        \ColA{0} & \ColB{1} & \ColB{0} & \ColB{0} & \ColA{0} & \ColB{1} & \ColB{0}
            & \ColB{0} & \ColA{0} & \ColA{0} \\
        \ColA{0} & \ColB{0} & \ColB{1} & \ColB{1} & \ColA{0} & \ColB{0} & \ColB{1}
            & \ColB{1} & \ColA{0} & \ColA{0} \\
        \ColA{0} & 0 & 0 & 0 & \ColA{0} & 0 & 0 & 0 & \ColA{0} & \ColA{1}
    \end{bmatrix}}.
\end{equation}
\end{footnotesize}%
Let us also endow $\LeftPAs$ with the right action $\SymmetricAction$ of the symmetric
groups such that, for any $u \in \SetPackedWords[n]$ and $\sigma \in \SetPermutations_n$,
\begin{equation}\label{Example_right_action}
    \BasisE_u \SymmetricAction \sigma
    = \BasisE_{\sigma^{-1}(u(1)) \; \dots \; \sigma^{-1}(u(\Length(u)))}.
\end{equation}
For instance,
\begin{math}
    \BasisE_{1411232} \SymmetricAction 3142 = \BasisE_{2322414}.
\end{math}
\medbreak

\begin{Theorem} \label{thm:left_packed_word_associative_symmetric_operad}
    The graded space $\LeftPAs$, endowed with the partial composition maps $\circ_i$ and the
    action $\SymmetricAction$ of the symmetric groups, is a symmetric operad.
\end{Theorem}
\medbreak

We call $\LeftPAs$ the \Def{left associative operad of packed words}. Observe that this
operad is not combinatorial since for any $n \geq 1$, there are infinitely many packed words
with $n$ as maximal value. However, as we shall see in the next sections, this operad
contains interesting quotient operads which are combinatorial.
\medbreak

%%%%%%%%%%%%%%%%%%%%%%%%%%%%%%%%%%%%%%%%%%%%%%%%%%%%%%%%%%%%%%%%%%%%%%%%%%%%%%%%%%%%%%%%%%%%
%%%%%%%%%%%%%%%%%%%%%%%%%%%%%%%%%%%%%%%%%%%%%%%%%%%%%%%%%%%%%%%%%%%%%%%%%%%%%%%%%%%%%%%%%%%%
%%%%%%%%%%%%%%%%%%%%%%%%%%%%%%%%%%%%%%%%%%%%%%%%%%%%%%%%%%%%%%%%%%%%%%%%%%%%%%%%%%%%%%%%%%%%
\section{Operadic quotients} \label{sec:quotients}
In this section, we regard the operad $\LeftPAs$ as set-operads through their
$\BasisE$-basis. This is possible since the composition of two elements of the
$\BasisE$-basis produces and element of the $\BasisE$-basis. For this reason, we shall write
here $u$ instead of $\BasisE_u$ for any $u \in \SetPackedWords$.
\medbreak

%%%%%%%%%%%%%%%%%%%%%%%%%%%%%%%%%%%%%%%%%%%%%%%%%%%%%%%%%%%%%%%%%%%%%%%%%%%%%%%%%%%%%%%%%%%%
%%%%%%%%%%%%%%%%%%%%%%%%%%%%%%%%%%%%%%%%%%%%%%%%%%%%%%%%%%%%%%%%%%%%%%%%%%%%%%%%%%%%%%%%%%%%
\subsection{Permutative congruences}
Let $P$ be a predicate on $\PositiveN^* \times \PositiveN \times \PositiveN \times
\PositiveN^*$. From $P$, we define the binary relation $\AdjacencyRelation$ on
$\PositiveN^*$ satisfying $uabv \AdjacencyRelation ubav$ for any $u, v \in \PositiveN^*$ and
$a, b \in \PositiveN$ such that $P(u, a, b, v)$ holds. Let also $\Equiv_P$ be the reflexive,
symmetric, and transitive closure of $\AdjacencyRelation$. The predicate $P$ is
\begin{enumerate}
    \item \Def{compatible with relabeling} if for any $u, v \in \PositiveN^*$ and $a, b \in
    \PositiveN$, $P(u, a, b, v)$ implies that $P\Par{f(u), f(a), f(b), f(v)}$ where $f :
    \PositiveN \to \PositiveN$ is any a strictly monotone map;
    \item \Def{compatible with subwords} if for any $u, v \in \PositiveN^*$ and $a, b \in
    \PositiveN$,  $P(u, a, b, v)$ implies that $P\Par{u', a, b, v'}$ for any $u', v' \in
    \PositiveN^*$ such that $u$ is a subword of $u'$ and $v$ is a subword of~$v'$.
\end{enumerate}
Observe that when $P$ is compatible with subwords, $\Equiv_P$ is a monoid congruence of the
free monoid on $\PositiveN$. Here are some examples of predicates compatible with relabeling
and subwords:
\begin{itemize}
    \item Let $\PredicateComm$ be the \Def{commutative predicate}, satisfying
    $\PredicateComm(u, a, b, v)$ for any $u, v \in \PositiveN^*$ and~\text{$a, b \in
    \PositiveN$.}
    \item Let $\PredicateSylv$ be the \Def{sylvester predicate}, satisfying
    $\PredicateSylv(u, a, c, v)$ for any $u, v \in \PositiveN^*$ and $a, c \in \PositiveN$
    such that there exists a letter $b$ in $v$ such that $a \leq b < c$. The equivalence
    relation $\Equiv_{\PredicateSylv}$ is the \Def{sylvester congruence} and has been
    introduced in~\cite{HNT05} to provide an alternative construction of the Loday-Ronco
    Hopf algebra.
    \item Let $\PredicateHypoplactic$ be the \Def{hypoplactic predicate}, satisfying
    $\PredicateHypoplactic(u, a, c, v)$ for any $u, v \in \PositiveN^*$ and $a, c \in
    \PositiveN$ such that there exists a letter $b$ in $u$ such that $a < b \leq c$ or there
    exists a letter $b$ in $v$  uch that $a \leq b < c$. The equivalence relation
    $\Equiv_{\PredicateHypoplactic}$ is the \Def{hypoplactic congruence} and has been
    introduced in~\cite{KT97}. This equivalence relation can be used to provide a
    construction of the Hopf algebra $\NCSF$ of noncommutative symmetric functions.
    \item Let $\PredicateBaxter$ be the \Def{Baxter predicate}, satisfying
    $\PredicateBaxter(u, a, c, v)$ for any $u, v \in \PositiveN^*$ and $a, c \in \PositiveN$
    such that there exist a letter $b$ in $u$ and a letter $b'$ in $v$ such that  $a \leq b'
    < b \leq d$ or $a < b \leq b' < d$. The equivalence relation $\Equiv_{\PredicateBaxter}$
    is the \Def{Baxter congruence} and has been introduced in~\cite{Gir12} to construct the
    Baxter Hopf algebra which contains the Loday-Ronco Hopf algebra as quotient.
\end{itemize}
In contrast, the \Def{plactic predicate} $\PredicatePlactic$, satisfying
$\PredicatePlactic(u, a, c, v)$ for any $u, v \in \PositiveN^*$ and $a, c \in \PositiveN$
such that $v$ is nonempty and its first letter $b$ is such that $a \leq b < c$ or $u$ is
nonempty and its last letter $b$ is such that $a < b \leq c$, is not compatible  with
subwords. The equivalence relation $\Equiv_{\PredicatePlactic}$ enjoys a lot of properties
(see for instance~\cite{DHT02} for properties related to the construction of Hopf algebras)
but does not play any role in this operadic context.
\medbreak

\begin{Theorem} \label{thm:left_packed_word_operad_permutative_congruences}
    If $P$ is a predicate compatible with relabeling and subwords, then $\Equiv_P$ is a
    nonsymmetric operad congruence of $\LeftPAs$.
\end{Theorem}
\medbreak

When $P$ satisfies the prerequisites of
Theorem~\ref{thm:left_packed_word_operad_permutative_congruences}, $\Equiv_P$ is a
\Def{permutative congruence}. Let us denote by $\theta_P : \LeftPAs \to \LeftPAs$ the map
sending any $u \in \LeftPAs$ to the minimal word w.r.t.\ the lexicographic order of the
$\Equiv_P$-equivalence class of $u$.
\medbreak

%%%%%%%%%%%%%%%%%%%%%%%%%%%%%%%%%%%%%%%%%%%%%%%%%%%%%%%%%%%%%%%%%%%%%%%%%%%%%%%%%%%%%%%%%%%%
%%%%%%%%%%%%%%%%%%%%%%%%%%%%%%%%%%%%%%%%%%%%%%%%%%%%%%%%%%%%%%%%%%%%%%%%%%%%%%%%%%%%%%%%%%%%
\subsection{Quotients forming combinatorial operads}
For any $k \geq 1$, let $\First{k} : \LeftPAs \to \LeftPAs$ be the map sending any $u \in
\LeftPAs$ to the packed word obtained by deleting any occurrence of a letter $a$ provided
that there are at least $k$ occurrence of $a$ on its left. For instance,
\begin{math}
    \First{2}(421431 \ColB{4} \ColB{1} 2) = 4214312.
\end{math}
Let us denote by $\Equiv_k$ the equivalence relation on $\LeftPAs$ satisfying, for any $u, v
\in \LeftPAs$, $u \Equiv_k v$ if $\First{k}(u) = \First{k}(v)$.
\medbreak

\begin{Proposition} \label{thm:left_packed_word_operad_first_congruences}
    For any $k \geq 1$, the equivalence relation $\Equiv_k$ is a nonsymmetric operad
    congruence of $\LeftPAs$.
\end{Proposition}
\medbreak

Since for any $k \geq 1$ and $n \geq 1$, there are finitely packed words having $n$ as
maximal value and where each letter appears at most $k$ times, the operad $\LeftPAs
/_{\Equiv_k}$ is combinatorial. The quotient $\LeftPAs /_{\Equiv_1}$ is the left associative
operad of permutations. The sequence of dimensions of $\LeftPAs /_{\Equiv_2}$ begins with
$2$, $14$, $222$, $6384$, $291720$, $19445040$ and forms Sequence~\OEIS{A105749}
of~\cite{oeis}.
\medbreak

A predicate $P$ is \Def{right-trivial} if for any $u, v \in \PositiveN^*$ and $a, b \in
\PositiveN$, $P(u, a, b, v)$ implies that $P(u, a, b, \epsilon)$.
\medbreak

\begin{Theorem} \label{thm:commutation_projection_permutative_first_k}
    If $P$ is a right-trivial predicate compatible with relabeling and subwords, then the
    maps $\First{k}$ and $\theta_P$ commute for any $k \geq 1$.
\end{Theorem}
\medbreak

Let us denote by $\Equiv_{P, k}$ the equivalence relation on $\LeftPAs$ satisfying, for any
$u, v \in \LeftPAs$, $u \Equiv_{P, k} v$ if $\theta_P\Par{\First{k}(u)} = \theta_P
\Par{\First{k}(v)}$. By Theorem~\ref{thm:commutation_projection_permutative_first_k},
$\Equiv_{P, k}$ is a nonsymmetric operad congruence of $\LeftPAs$. Let us review some
quotients of $\LeftPAs$ obtained via such congruences:
\begin{itemize}
    \item The commutative predicate $\PredicateComm$ is right-trivial. The nonsymmetric
    operad $\LeftPAs /_{\Equiv_{\PredicateComm, 1}}$ is the nonsymmetric associative operad,
    and for any $k \geq 1$ and any $n \geq 1$, the dimension of $\LeftPAs
    /_{\Equiv_{\PredicateComm, k}}(n)$ is $k^n$. A set of representatives of this quotient
    is the set of weakly increasing packed words having $n$ as maximal value and such that
    each maximal factor of equal letters is of length at most $k$.
    \item The sylvester predicate is not right-trivial, but by setting $\PredicateCoSylv$ as
    the predicate satisfying $\PredicateCoSylv(u, a, b, v)$ if and only if
    $\PredicateCoSylv(v, a, b, u)$, $\PredicateCoSylv$ is right-trivial. The sequence of the
    dimensions of $\LeftPAs /_{\Equiv_{\PredicateCoSylv, 1}}$ begins with $1$, $2$, $5$,
    $14$, $132$, $429$ and forms Sequence~\OEIS{A000108} of~\cite{oeis} of Catalan numbers.
    \medbreak

    \begin{Proposition}
        The nonsymmetric operad $\LeftPAs /_{\Equiv_{\PredicateCoSylv, 1}}$ is isomorphic to
        the duplicial operad.
    \end{Proposition}
    \medbreak

    The sequence of the dimensions of $\LeftPAs /_{\Equiv_{P, 2}}$ begins with $2$, $10$,
    $66$, $498$,  $4066$, $34970$ and forms Sequence~\OEIS{A027307} of~\cite{oeis}. The
    operads $\LeftPAs /_{\Equiv_{\PredicateCoSylv, k}}$ are generalizations of the duplicial
    operad.
\end{itemize}
\medbreak

%%%%%%%%%%%%%%%%%%%%%%%%%%%%%%%%%%%%%%%%%%%%%%%%%%%%%%%%%%%%%%%%%%%%%%%%%%%%%%%%%%%%%%%%%%%%
%%%%%%%%%%%%%%%%%%%%%%%%%%%%%%%%%%%%%%%%%%%%%%%%%%%%%%%%%%%%%%%%%%%%%%%%%%%%%%%%%%%%%%%%%%%%
%%%%%%%%%%%%%%%%%%%%%%%%%%%%%%%%%%%%%%%%%%%%%%%%%%%%%%%%%%%%%%%%%%%%%%%%%%%%%%%%%%%%%%%%%%%%
\section{Analogue of Dynkin idempotents to Zie} \label{sec:dynkin}
The weak  order on permutations has several analogues for packed words. The \Def{right weak
order} for packed words is defined via the following covering relation: for $u,v \in
\SetPackedWords_k$ $u \prec_r v$ if and only if $v = u \cdot \tau$, for some transposition
$\tau \in \SetPermutations_k$ and $\Inversions(v) = \Inversions(u)+1$. Given $u \in
\SetPackedWords$, the \Def{set of left-inversions} $\text{Inv}_\ell(u)$ of $u$ is the set of
pairs $(i,j)$ such that $i<j$ and all appearances of $i$ in $u$ occur after all appearances
of $j$ in $u$. The \Def{left weak order} of packed words is defined via the following
covering relation: for $u,v \in \SetPackedWords[n]$, $u \prec_\ell v$ if and only if $v
=\tau \cdot u$, for some transposition $\tau \in \SetPermutations_n$ and
$|\text{Inv}_\ell(v)| =  |\text{Inv}_\ell(u)|+1$. 
\medbreak

Consider the following graded bases for $\RightPAs$:
\begin{equation}
    \BasisE_u = \sum_{v \leq_r u}\BasisF_v
    \quad \text{and} \quad 
    \BasisF_u = \sum_{u \leq_\ell v} \BasisM_v.
\end{equation}
The set $\{\BasisM_u : u \in \SetPackedWords_n, n \geq 0\}$ is a linear basis of the
primitive space of the Hopf algebra $\WQSym$ of packed words (see, for example,
\cite{BGLPV2023}). Following the terminology of \cite{AM2017}, we called this space the
\Def{space of Zie elements}. In the case of  permutations, a \Def{Lie element} is an element
of the primitive space of $\FQSym$, the Hopf algebra of permutations~\cite{MR1995,DHT02}.
\medbreak

It is straightforward to show that the $\BasisF$-basis satisfies the same internal operation
in $\RightPAs$ of the $\BasisE$-basis, as described in~\eqref{RightPAs_Def}. Given $u, v \in
\SetPackedWords$, define a new operation on $\RightPAs$ as
\begin{equation}
    \left\{\BasisF_u,  \BasisF_v\right\}
    := (\BasisM_{12} \circ_2 \BasisF_v)\circ_1 \BasisF_u.
\end{equation}
For instance,
\begin{math}
    \{\BasisF_{11}, \BasisF_{121}\}
    = ((\BasisF_{12}-\BasisF_{21})\circ_2 
    \BasisF_{121})\circ_1 \BasisF_{11} = \BasisF_{11232} - \BasisF_{33121}.
\end{math}
Remark that this operation is not the induced Lie bracket of $\WQSym$; while $\{\BasisF_1,
\BasisF_1\}= \BasisF_{12} - \BasisF_{21}$, we have $[\BasisF_1, \BasisF_1]=0$.
\medbreak

Now, given $k \geq 1$, the \Def{Dynkin idempotent} $\frac{1}{k}\, \theta_k$ is defined as
the unique element in $\RightAs(k)$ for which $\theta_k \cdot a_1a_2a_3 \cdots a_k = [\cdots
[[a_1,a_2], a_3], \cdots, a_k]$ is the $k-1$ left nested commutator bracket in the  vector
space $\mathbb{K}\langle \mathbb{N} \rangle$, for any word $a_1a_2a_3 \cdots a_k \in
\mathbb{K}\langle \mathbb{N} \rangle$. Here, we consider the natural right action of
permutations (as in \eqref{Example_right_action}). It is shown in \cite[Lem. 5.2]{AL2007}
that we can express $\theta_k$ as 
\begin{math}
    \theta_k = \{ \dots \{ \{\BasisF_1, \BasisF_1\}, \BasisF_1 \}, \dots , \BasisF_1\}
\end{math}
(there are $k-1$ left nested brackets). Given a composition $\mathsf{c}=c_1c_2\cdots c_n \in
\mathbb{N}^*$ of sum $k$, we define the \Def{$c$-Dynkin's idempotent}
$\frac{1}{k}\theta^{\mathsf{c}}_k$, where 
\begin{equation}
    \theta_k^{\mathsf{c}}:= \left\{ \dots \left\{ \left\{\BasisF_{\text{id}_{c_1}},
    \BasisF_{\text{id}_{c_2}}\right\}, \BasisF_{\text{id}_{c_3}} \right\}, \dots,
    \BasisF_{\text{id}_{c_n}}\right\}
\end{equation}
and $\text{id}_r$ denotes the packed word $11\cdots 1$ of length $r$. More generally, let
$\text{id}(\mathsf{c}) \in \SetPackedWords_k[n]$ be the non-decreasing packed word with
$c_r$ copies of $r$, for every $r \in [n]$.
\begin{Theorem}
    For any $k \geq 1$,
    \begin{equation}
        \theta_k^{\mathsf{c}}
        = \sum_{\substack{
            \text{id}(\mathsf{c}) \; \leq_\ell \; u \\
            u = \text{id}\left(c_1\right) 2 \; \circ_2 \; u'}} 
        \BasisM_u.
    \end{equation}
    In particular, $\theta^{\mathsf{c}}_k$ is a Zie element in $\WQSym$.
\end{Theorem}
\medbreak

%%%%%%%%%%%%%%%%%%%%%%%%%%%%%%%%%%%%%%%%%%%%%%%%%%%%%%%%%%%%%%%%%%%%%%%%%%%%%%%%%%%%%%%%%%%%
%%%%%%%%%%%%%%%%%%%%%%%%%%%%%%%%%%%%%%%%%%%%%%%%%%%%%%%%%%%%%%%%%%%%%%%%%%%%%%%%%%%%%%%%%%%%
%%%%%%%%%%%%%%%%%%%%%%%%%%%%%%%%%%%%%%%%%%%%%%%%%%%%%%%%%%%%%%%%%%%%%%%%%%%%%%%%%%%%%%%%%%%%
\section{Open questions and future work} \label{sec:open_questions}
Here are some questions and possible ways to continue this work.
\begin{enumerate}[label={\bf (\arabic*)}]
    \item \textbf{Bases of $\RightPAs$ and partial orders on packed words} --- One of the
    initial motivations of this work was to mimic the constructions relating the associative
    operad and the weak order developed by Aguiar and Livernet in~\cite{AL2007}. They showed
    that the partial composition on the $\BasisM$-basis of $\As$ is encoded by intervals,
    which is not the case for $\RightPAs$. This leads to the exploration of other posets
    structures on packed words.
    \item \textbf{Quotients of $\RightPAs$} --- A second motivation of this work was to
    construct at the level of our operads of packed words similar structures as the
    duplicial and interstice operads, which are quotients of $\As$~\cite{AL2007}. This has
    been done for $\LeftPAs$. The question remains open for $\RightPAs$. We conjecture that
    the analog of the duplicial operad, as quotient of $\RightPAs$, involves Schröder trees.
    There are other operads on Schröder trees defined in~\cite{Chapoton2002}
    and~\cite{Gir2015}. A question is to investigate whether these operads are isomorphic.
    \item \textbf{Generalizations of $\RightPAs$} --- The partial compositions of the
    $\RightPAs$ can be extended on generalizations of packed words, while preserving the
    structure of a nonsymmetric operad. Such extensions work both for the sets of parking
    functions and of endofunctions. The study of these two new structures is planned. In
    contrast, the similar extensions for $\LeftPAs$ do not form operads.
    \item \textbf{Quotients of $\LeftPAs$} --- The study of the quotients $\LeftPAs
    /_{\Equiv_{\PredicateComm, k}}$ and $\LeftPAs /_{\Equiv_{\PredicateCoSylv, k}}$ of
    $\LeftPAs$ is planed. This last quotient involves Schröder trees subjected to some
    conditions and with a unusual notion of arity. This further study of operads encompasses
    the description of presentations by generators and relations and combinatorial
    descriptions of their partial compositions on right combinatorial objects.
    \item \textbf{Diagonal of the permutahedron} --- Given a polytope $\mathtt{P}$, it is a
    general problem to find a cellular approximation which is homotopic to the diagonal map
    $\Delta_{\mathtt{P}}: \mathtt{P} \to \mathtt{P} \times \mathtt{P}$, and which agrees
    with $\Delta_{\mathtt{P}}$ on the vertices of $\mathtt{P}$. For instance, the cellular
    chains on the associahedra are  endowed with an operad structure which encodes
    associative algebras up to homotopy. A natural question is how to relate our two operad
    structures $\LeftPAs$ and $\RightPAs$ to coherent cellular approximations of the
    diagonal of the permutahedra.
\end{enumerate}
\medbreak

%%%%%%%%%%%%%%%%%%%%%%%%%%%%%%%%%%%%%%%%%%%%%%%%%%%%%%%%%%%%%%%%%%%%%%%%%%%%%%%%%%%%%%%%%%%%
%%%%%%%%%%%%%%%%%%%%%%%%%%%%%%%%%%%%%%%%%%%%%%%%%%%%%%%%%%%%%%%%%%%%%%%%%%%%%%%%%%%%%%%%%%%%
%%%%%%%%%%%%%%%%%%%%%%%%%%%%%%%%%%%%%%%%%%%%%%%%%%%%%%%%%%%%%%%%%%%%%%%%%%%%%%%%%%%%%%%%%%%%
\printbibliography

\end{document}